\input amstex
\documentstyle{amsppt}
\vsize=7.0in \hsize 6.8 in
\loadbold \topmatter
\title On singular connections and geometrically atomic maps\endtitle
\author Valery Marenich and Karel Packalen \endauthor
\thanks
\endthanks
\address Department of Tehnology, The University of Kalmar, 905, S-391~82, Kalmar, Sweden.\endaddress
\email valery.marenich\@ hik.se, karel.packalen\@ hik.se \endemail \abstract
We prove that an introduction of the class of geometrically atomic bundle maps by Harvey and Lawson in their theory of singular
connections is not necessary because an arbitrary map satisfies the conditions of geometric atomicity.
\endabstract
\subjclass  Primary 53C07, 49Q15; Secondary 58H99, 57R20. Research was supported by {Vetenskapr{\aa}det} (Sweden) and the
Faculty of Natural Sciences of the University of Kalmar, (Sweden)
\endsubjclass \keywords singular connections, geometrically atomic maps \endkeywords
\endtopmatter

\document
\baselineskip 12pt

\head \bf 1. Singular connections \endhead

Let us recall some basic facts about the Harvey-Lawson theory of singular connections, see [HL1]\footnote{See also [HL2]-[HL4],
[N] and [Z]}. We consider a smooth bundle map $\alpha: E\to F$ between complex vector bundles with connections $D^E$ and $D^F$
correspondingly over the same manifold $X$. We denote by $\nu=dim(X)$, for a point $x$ of $X$ the fibers over $x$ are $E_x$ and
$F_x$, $\alpha_x: E_x\to F_x$ is the restriction of $\alpha$ to these fibers, while $m$ and $n$ denote complex dimensions of
these fibers.

For the Chern-Weil characteristic forms $\Phi(\Omega)$ of $(E,D^E)$ and $(F,D^F)$ and the class of "geometrically atomic" maps
$\alpha$ the authors in [HL1] establish a canonical co-homology
$$
\Phi(\Omega) - \sum_{k\geq 0} Res_{\Phi,k}[\Sigma_k(\alpha)] = dT, \tag 1.1
$$
where $\Sigma_k(\alpha)=\{x\in X : dim ker(\alpha)=k\}$, $Res_{\Phi,k}$ is a smooth residue form along $\Sigma_k(\alpha)$, and
$T$ is a canonical $L^1_{loc}$-form on $X$. When $m=n$ the last formula can be written
$$
\Phi(\Omega^F) - \Phi(\Omega^E)= \sum_{k>0} Res_{\Phi,k}[\Sigma_k(\alpha)] = dT. \tag 1.2
$$

A bundle morphism $\alpha: E\to F$ can be considered as a cross-section of the vector bundle $\pi: Hom(E,F)\to X$. Denote by
$\beta: \pi^*E\to \pi^*F$ a tautological morphism between bundles $ \pi^*E$ and $\pi^*F$ over the total space of the bundle
$Hom(E,F)$ (which at $A\in Hom(E,F)$ is given by itself). There is a natural compactification $Hom(E,F)\subset G$ where $\pi:
G\equiv G_m(E\oplus F)\to X$ is the Grassmann bundle of complex m-planes in $E\oplus F$ while we assign to $A: E_x\to F_x$ its
graph $V_A$ in $E_x\oplus F_x$.

The key object in the construction of the co-cycle $T$ is the following {\bf multiplicative flow} $\phi_t: E\oplus F \to
E\oplus F$ defined by $\phi_t(e,f)=(te,f)$ for $0<t<\infty$. If $P^1(R)=R\cup\{\infty\}$ is the projective line, $G^{\oplus 2}$
the fibre product of $G$ with itself over $X$ then
$$
{\Cal T}=\{(t,\phi_t(V),V)\in P^1(R)\times G^{\oplus 2} | 0<t<\infty \text{ and } V\in G\} \tag 1.3
$$
is called the total graph of the flow. Denote by [${\Cal T}$] the current given by integration over ${\Cal T}$ and define
$T=pr_*[{\Cal T}]$ where $pr:P^1(R)\times G^{\oplus 2}\to G^{\oplus 2} $ is the projection.

As the authors wrote in [HL4] the main observation they had made was the following statement.

\proclaim{Theorem~1. (see Proposition 3.1 in [HL4])} ${\Cal T}$ is a sub-manifold of finite volume in $P^1(R)\times G^{\oplus
2}$ over each compact subset of $X$.\footnote{It was concluded from the real analyticity of the multiplicative flow with the
help of general results from a geometric measure theory, see [F].}
\endproclaim

\medskip

The theory of singular connections was developed for a class of so called geometrically atomic bundle maps. The bundle map
$\alpha: E\to F$ is called {\bf geometrically atomic} if the subset
$$
T_\alpha = \{({{1}\over{t}}\alpha_x,\alpha_x) \in G^{\oplus 2}| x\in X, 0<t<\infty\} \tag 1.4
$$
has locally finite $(\nu+1)$-dimensional measure in $G^{\oplus 2}$ (recall that $\nu=dim(X)$).

The set $T_\alpha$ can be divided into two sets: the set $T^s_\alpha$ consisting of points where $\alpha$ has zeros and
therefore having locally finite $\nu$-dimensional measure; and the set $T^r_\alpha$ which is a line bundle over the set
$X\backslash \text{Zero}(\alpha)$. This last set is an open subset of $X$ while $T^r_\alpha$ is a smooth sub-manifold which
fibers over $X\backslash \text{Zero}(\alpha)$ with a fiber equals the trajectory $({{1}\over{t}}\alpha_x,\alpha_x)$ of
$\alpha_x$ under the multiplicative flow. Therefore, by Fubini theorem the sub-manifold $T^r_\alpha$ has locally finite
$(\nu+1)$-dimensional measure if all these fibers $\{\phi_t(\alpha_x), 0<t<\infty\}$ have uniformly bounded length considered
as curves in Grassmann manifolds $G^{\oplus 2}_x$. In this note we prove this uniform boundedness of lengths of such curves in
Grassmannians, see Theorem~3 below. This will imply our main result as follows.

\proclaim{Theorem~2} All bundle maps $\alpha: E\to F$ are geometrically atomic.
\endproclaim

\head \bf 2. Linear curves of subspaces in Grassmannians\endhead

Denote $p=2m$, $q=2n$. Then every $G^{\oplus 2}_x$ is $G_{p+q}^p$ - the Grassmann manifold of all $p$-dimensional subspaces of
the $(p+q)$-dimensional Euclidean space. We consider $G_{p+q}^p$ with a natural Riemannian metric $g$ such that the
factorization $\Pi : {O}(p+q)\to {O}(p+q)/{O}(p)\times {O}(q)=G_{p+q}^p$, where the orthogonal group $O(p+q)$ is provided with
the Lipcshitz-Killing metric, is a Riemannian submersion.

We say that a curve $V(t), 0\leq t \leq T$ in $G_{p+q}^p$ of $p$-subspaces is {\bf linear} if the subspaces $V(t)$ are
generated by some vectors $e_i(t), i=1,...,k$ linearly depending on $t$, i.e., such that
$$
e_i(t)=e_i+ t d_i. \tag 2.1
$$
In our case of the trajectories of the multiplicative flow are the following linear curves:
$$
V(t)=\phi_t(V_A)=({{1}\over{t}}\alpha_x,\alpha_x)
$$
where the $p$-subspace $V_A$ is the graph of the map $A=\alpha_x$. If $\{e_i,1\leq i\leq p\}$ is some base of the first factor
$R^p$ of $R^{p+q}=R^p\times R^q$ (where $p=2m$, $q=2n$) then $V(1)=V_A$ and is generated by $e_i+Ae_i$, while $V(t)$ is
generated by the vectors
$$
e_i(t)=e_i+ t Ae_i. \tag 2.2
$$
Denote $d_i=Ae_i$ and call these vectors derivatives.

As we said above our main technical result is:

\proclaim{Theorem~3} The lengths of all linear curves $V(t)$ in a Grassmann manifold $G^p_{p+q}$ are uniformly bounded.
\endproclaim

\demo{Proof} First we verify that geodesic curvatures of linear curves are uniformly bounded. Consider the point $V(t_0)$ on
the linear curve $V(t)$ and estimate its geodesic curvature at this point. We split the space $R^{p+q}$ as $R^{p+q}=R^p\times
R^q$ where the first factor $R^p$ of $R^{p+q}=R^p\times R^q$ coincide with $V(t_0)$ and will be called the {\bf vertical
subspace} and its vectors - vertical, the second factor $R^q$ of $R^{p+q}=R^p\times R^q$ which is orthogonal to the first one
we call the {\bf horizontal} subspace and its vectors - horizontal vectors. Vertical and horizontal components of the vector
$v$ will be denoted below as ${\Cal V}v$ and ${\Cal H}v$ correspondingly.

Changing $t$ to $(t-t_0)$ if necessary in (2.2) we assume  from now on that $t_0=0$ and denote $V(0)=V$. In return the
derivatives $d_i$ in a representation
$$
V(t)=\{e_i(t)=e_i+td_i | 1\leq i\leq p \} \tag 2.2
$$
may be no longer horizontal vectors as before. Continue the correspondence $e_i\to {\Cal V}d_i$ to the linear map ${\Cal
V}D:V\to V$ and take its polar representation ${\Cal V}D=BU$ where $U:V\to V$ is some orthogonal map $U=(u_i^s)$ from $O(p)$,
and $B$ a diagonal in some orthonormal base.\footnote{$B$ can be assumed to be positive, see the construction in the Lemma~3
below of the partition $\{-\infty<t_1<...<t_m<\infty\}$.} The pre-image of this base under transformation $U$ we denote by
$\{e_1, ... , e_p\}$ of $V$. Complete this base to some orthonormal base $\{e_1, ... , e_p, e_{p+1}, ..., e_{p+q}\}$ of
$R^{p+q}$ and continue the linear map $D:V\to R^{p+q}$ sending $e_i$ to $d_i$ to the map $D:R^{p+q}\to R^{p+q}$ given by a
matrix
$$
D=\left(
\matrix B           & ({\Cal H}D)\\
       -({\Cal H}D)^*          & I\\
\endmatrix \right)
\left(
\matrix U  & 0 \\
                 0              & I\\
\endmatrix \right) . \tag 2.3
$$

Now we see that the linear curve $V(t)$ is the image of the subspace generated by the first $p$ coordinate vectors of the
orthonormal base $\{e_1, ... , e_p, e_{p+1}, ..., e_{p+q}\}$  in $R^{p+q}$ under the map with the following matrix:
$$
\left(
\matrix I+t({\Cal  V}D) & t({\Cal H}D)\\
                 -t({\Cal H}D)^*  & I\\
\endmatrix \right)
\left(
\matrix U  & 0 \\
                 0              & I\\
\endmatrix \right) . \tag 2.4
$$
where ${\Cal V}D=BU$ and ${\Cal H}D$ are the matrices of vertical and horizontal components of the derivatives $d_i$
correspondingly: ${\Cal V}D=(d_i)^j, i=1, ... , p; \quad  j=1,...,p$ and ${\Cal H}D=(d_i)^j, i=1, ... , p; \quad
j=p+1,...,p+q$.

Next we define another and "nearly" orthogonal transformation of $R^{p+q}$ mapping $V$ into $V(t)$. First we divide every
vector $e_i, 1\leq i\leq p$ by its norm after $I+tBU$, i.e., by $(1+t\lambda_i)$. We denote this map
$C(t)=diag((1+t\lambda_1)^{-1},...(1+t\lambda_p)^{-1})$. Its composition with the map given by (2.4)
$$
\left(
\matrix I+t({\Cal  V}D) & t({\Cal H}D)\\
                 -t({\Cal H}D)^*  & I\\
\endmatrix \right)
\left(
\matrix U  & 0 \\
                 0              & I\\
\endmatrix \right)
\left(
\matrix C(t)  & 0 \\
                 0              & I\\
\endmatrix \right)  \tag 2.5
$$
also send $V$ onto $V(t)$. By a direct calculation it is easy to see that this last map is given in the orthonormal base
$\{e_i, 1\leq i\leq p+q\}$ by the matrix
$$
A(t)=\left(
\matrix \tilde U & t\tilde D\\
                 -t\tilde D^*  & I\\
\endmatrix \right)
\tag 2.6
$$
where $\tilde U$ is some orthogonal matrix and for the components of $\tilde D=(\tilde d_i^j), 1\leq i\leq p, p+1\leq j\leq
p+q$ it holds that
$$
\tilde d_i^j=(1+t\lambda_i)^{-1} u_{i}^s d_s^j. \tag 2.7
$$
Another calculation shows that the vectors $\tilde e_i(t)=A(t)e_i, 1\leq i\leq p+q$ are "nearly" orthonormal:
$$
|(\tilde e_i(t),\tilde e_j(t))-\delta_{ij}|\leq t^2\|\tilde D\|^2, \tag 2.8
$$
where
$$
\|\tilde D\|^2=\sum\limits_{i=1,...,k;j=k+1,...,n}((\tilde d_i)^j)^2. \tag 2.9
$$
For small $t$ it holds $(1+t\lambda_i)^{-1}\leq 1+t\Lambda_1$ for $\Lambda_1=2 max\{-\lambda_i|1\leq i\leq p\}$. Therefore,
from (2.7) we see
$$
\|\tilde D\|^2 \leq (1+t\Lambda_1)\|{\Cal H} D\|^2, \tag 2.10
$$
or
$$
|(\tilde e_i(t),\tilde e_j(t))-\delta_{ij}|\leq t^2(1+t\Lambda_1)\|{\Cal H} D\|^2, \tag 2.11
$$
where
$$
\|{\Cal H} D\|^2=\sum\limits_{i=1,...,k;j=k+1,...,n}((d_i)^j)^2. \tag 2.12
$$
Now we apply inductively to $\{\tilde e_i(t), i=1,...p+q\}$ an orthogonalization process as follows. Denote $W_1(t)$ the
subspace generated by all $\{\tilde e_i(t), i<p+q\}$ and find $\bar e_{p+q}(t)$ unit and normal to it. From (2.11) it follows
that
$$
\|\tilde e_{p+q}(t)-\bar e_{p+q}(t)\| \leq t^2(1+t\Lambda^1_1)\|{\Cal H} D\|^2 \tag 2.13
$$
for some $\Lambda^1_1$ depending on $\Lambda_1$. Next find $\bar e_{p+q-1}(t)$ in $W_1(t)$ which is unit and normal to the
subspace $W_2(t)$ generated by all $\{\tilde e_i(t), i<p+q-1\}$. Again from (2.11) and (2.13) we get
$$
\|\tilde e_{p+q-1}(t)-\bar e_{p+q-1}(t)\| \leq t^2(1+t\Lambda^2_1)\|{\Cal H} D\|^2 \tag 2.14
$$
for some $\Lambda^2_1$ depending on $\Lambda_1$ and $\Lambda^1_1$. And so on. Finely we construct the orthonormal base $\{\bar
e_i(t), i=1,...p+q\}$ of $R^{p+q}$ such that
$$
\|\tilde e_i(t) - \bar e_i(t)\| \leq (1+t\Lambda_2) t^2\|{\Cal H}D\|^2 \tag 2.15
$$
where the constant $\Lambda_2$ again depends only on $\Lambda_1$ and the dimension $(p+q)$. By construction the first $p$
vectors of $\{\bar e_i(t), i=1,...n\}$ generate $V(t)$, or $V(t)$ is the image of $V(0)$ under the orthogonal transformation
with a matrix:
$$
O(t)=\left(
\matrix  \tilde U         & t({\Cal H}D)\\
                 -t({\Cal H}D)  & I\\
\endmatrix
\right) +t^2 G, \tag 2.16
$$
where for the norm of the matrix $G$ from (2.15) it follows
$$
\|G\|\leq(1+ t\Lambda_3) \|{\Cal H}D\|^2, \tag 2.17
$$
for some constant $\Lambda_3$ depending on $\Lambda_2$. Finally note that the curve $V(t)$ in $G^p_{p+q}$ is the image under
submersion $\Pi:O(p+q)\to G^p_{p+q}$ of the curve $O(t)$ of orthogonal transformations in $O(p+q)$ so that the vector $\dot
V(0)$ is the image of $\dot O(0)$ under the differential $\Pi_*$. From (2.14) we see that the tangent vector of the curve of
the orthogonal transformations $O(t)$ at the moment $t=0$ is
$$
\dot{O}(0)=\left(
\matrix 0 & {\Cal H}D\\
                 -{\Cal H}D  & 0\\
\endmatrix
\right) . \tag 2.18
$$
By definition of the Lipschitz-Killing metric the length of this vector equals $\|{\Cal H}D\|$. Note also that this vector is
the horizontal one for the Riemannian submersion $\Pi : {O}(p+q)\to {O}(p+q)/{O}(p)\times {O}(q)=G_{p+q}^p$, so that it has the
same length as its image $\dot V(0)$ under the differential $\Pi_*$. Therefore, the length of the vector $\dot V(0)$ in
$G^p_{p+q}$ equals the length of $\dot{O}(0)$ in the Lipschitz-Killing metric on ${O}(p+q)$, i.e., is exactly $\|{\Cal H}D\|$;
and we arrive at the following statement.

\proclaim{Lemma~1}
$$
\|\dot V(0)\|^2=\|{\Cal H}D\|^2=\sum\limits_{i=1,...,k;j=k+1,...,n}((d_i)^j)^2.
$$
\endproclaim

\medskip

\proclaim{Lemma~2} The geodesic curvature $K(V(t))$  of the linear curve $V(t)\subset G^p_{p+q}$ is uniformly bounded.
\endproclaim

\demo{Proof} To prove this we use a well-known equality between a geodesic curvature of some curve in a Riemannian manifold and
a distance from this curve to the geodesic, issuing from the same initial point with the same velocity vector. It is known that
the geodesic curvature $K(V(0))$ of the curve $V(t)$ and the norm of its vector of velocity are related by the following
equality:
$$
K(V(0))=  \lim_{t\to 0}2{{dist_{G^p_{p+q}}(V(t), \bar V(t))}\over{\|\dot V(0)\|^2 t^2}}, \tag 2.19
$$
where $\bar V(t)$ is the geodesic in $G^p_{p+q}$ issuing from $V(0)$ with the same velocity vector: $\dot{\bar V(0)}=\dot
V(0)$. Every geodesic  $\bar V(t)$ in the base $G^p_{p+q}$ of the Riemannian submersion $\Pi : {O}(p+q)\to
{O}(p+q)/{O}(p)\times {O}(q)=G^p_{p+q}$ is the image under the map $\Pi$ of some horizontal geodesic in ${O}(p+q)$, i.e., some
1-parameter subgroup of orthogonal transformations of $R^{p+q}$. Consider the 1-parameter subgroup $\bar O(t)$  of $O(p+q)$
generated by the vector $\dot O(0)$ above:
$$\dot O(0)=\left(
\matrix 0 & {\Cal H}D\\
      -{\Cal H}D  & 0\\
\endmatrix
\right).
$$
Because the vector $\dot O(0)$ is horizontal the geodesic $\bar O(t)$ in $O(p+q)$ is a horizontal geodesic and goes under
submersion $\Pi$ onto some geodesic in $G^p_{p+q}$. Because, as we verified above, the vector $\dot O(0)$ has the image under
the differential $\Pi_*$ which is equal to $\dot V(0)$, this geodesic has the same velocity vector at initial point $t=0$,
i.e., we conclude $\bar V(t)=\Pi(\bar O(t))$. A direct calculation shows that
$$
{\bar O}(t)=\exp(t\dot {O}(0))=I+t\dot {O}(0)+t^2\bar{G} + ... \tag 2.20
$$
where for some constant $\Lambda_4$ we have
$$
\|\bar {G}\|\leq \Lambda_4 \|{\Cal H}D\|^2. \tag 2.21
$$
Because of (2.14) and (2.15) above this gives the inequality
$$
\|\bar {O}(t) - {O}(t)\|\leq t^2\Lambda_5\|{\Cal H}D\|^2/2 \tag 2.22
$$
for some constant $\Lambda_5$ depending on $\Lambda_3$ and $\Lambda_4$. Because (as every Riemannian submersion) $\Pi$ does not
increase distances the last formula yields an inequality
$$
\|\bar {V}(t) - {V}(t)\|\leq t^2\Lambda_5\|{\Cal H}D\|^2/2 \tag 2.23
$$
implying the claim of the lemma. Lemma~2 is proved.
\enddemo

The next Lemma is the  main technical point of our arguments.

\proclaim{Lemma~3} Let $V(t), 0\leq t \leq T$ be some linear curve in $G^p_{p+q}$ such that its geodesic curvature $K(V(t))$ is
bounded by some constant $\Lambda_5$ for all $0\leq t\leq T$. Then the length of $V(t)$ is bounded by some constant $\Lambda_6$
which depends on $\Lambda_5$, but does not depend on $T$.
\endproclaim

\demo{Proof} Proof of the Lemma follows from  simple compactness and monotonicity arguments as follows.

In addition to the Riemannian metric $g$ on $G^p_{p+q}$ we used above (coming from the Lipschitz~-~Killing metric on ${O}(p+q)$
under the submersion $\Pi$), the Grassmanian $G^p_{p+q}$ also admits the following non-Riemannian metric $\angle$ ("angle"):
for two subspaces $V$ and $W$ in $R^{p+q}$ the angle between them is
$$
\angle(V,W) = \sup\limits_{v\in V}\inf\limits_{w\in W}{\angle(v,w)}. \tag 2.24
$$
Clearly, we always have $0\leq\angle(V,W)\leq \pi/2$ and $\angle(V,W)=\pi/2$ if and only if some vector $w$ of $W$ is
orthogonal to $V$.

Let $V(t)$ be our linear curve  determined by $\{e_i,d_i;i=1,...,p\}$ as above; i.e., $V(t)=(I+tA)V$ for some linear map
$A:V\to R^{p+q}$. For a vector $e(0)$ of $V(0)$ we have $(e(0),e(t))=0$ for $e(t)=(I+tA)e(0)$ from $V(t)$ if and only if
$$
(e(0), (I+tA)e(0))=0, \tag 2.25
$$
Note that there exists only a finite number of moments $t_i$ such that (2.25) holds for some non-zero $e(0)$ of $V(0)$. Indeed,
if we denote by $B:R^{p+q}\to V(0)$ the orthogonal projection onto $V(0)$ then such vectors belong to the kernel of the
operator $B(I+tA):V(0)\to V(0)$ which is non-zero if and only if the determinant function $P(t)=det(B(I+tA))$ vanishes. Because
this determinant function is a polynomial of degree $p$ and is not identically zero (note that $P(0)=1$) it has at most $p$
zeros $t_i, 1\leq i\leq m\leq p$.

For a partition $-\infty<t_1<...<t_m<\infty$ it is not difficult to see that for arbitrary $t_i < t'<t"< t_{i+1}$ we have
$$
(e(t'),e(t"))>0. \tag 2.26
$$
because otherwise (as is easy to check) for some $t^*$ from the interval $(t',t")$ it would follow $(e(0),e(t^*))=0$.

Condition (2.26) implies that the angle between some fixed vector $e(t')$ and $e(t")$ is a monotonely increasing function on
$t"$ when $t'<t"$ are from $(t_i,t_{i+1})$. The same condition implies, that the angle between  $e(t")$ and $V(t')$ is also a
monotonely increasing function on $t"$ under the same restriction. Indeed, by definition for all $t$ we have $e(t)=e(0)+ td$
for some $e(0)$ from $V(0)$. Let $d=v+w$, where $v$ belongs to the subspace $V(t')$, and $w$ is normal to it. Let also
$e(t')=a+b$, where the vector $a=\lambda v$ is parallel to $v$, and $b$ is normal to $v$. Then $e(t")=b + (\lambda + (t"-t'))v
+ (t"-t')w$, and (2.26) means that
$$
(e(t'),e(t')) + (t"-t')(e(t'),v)>0. \tag 2.27
$$
The component of $e(t")$ normal to $V(t')$ equals $(t"-t')w$  while the tangent component of $e(t")$ to $V(t')$ is
$e(t')+(t"-t')v$. So for the angle $\phi(t")$ between $e(t")$ and the subspace $V(t')$ we have
$$
tg(\phi(t"))={{\sqrt{(w,w)}(t"-t')}\over{\sqrt{(e(t')+(t"-t')v, e(t')+(t"-t')v)}}},
$$
and by a direct calculation we see, that due to (2.27) the derivative of $\phi(t")$ is positive. Because this is true for an
arbitrary $e(t")$ we conclude that the angle between $V(t")$ and $V(t')$ is also a monotonely increasing function on $t"$ when
$t'<t"$ are from $(t_i,t_{i+1})$. To prove this it is sufficient to note that the derivative of the angle between $V(t')$ and
$V(t")$  on $t"$ equals the maximum of derivatives of angles between  $V(t')$ and those vectors $e(t")$  of $V(t")$ which have
a maximum angle with $V(t')$.

Thus,  every linear curve $V(t)$ we can divide into not more than $p+1$ intervals $\{V(t)| t_i<t<t_{i+1}\}$ such that the
"angle"-function $\phi_s(t)=\angle(V(s), V(s+t))$ is monotonely increasing for $0<t<t_{i+1}-s$ and monotonely decreasing for
$t_i-s<t<0$. The claim of the Lemma will follow if we prove it for every such interval. In order to do this we note that this
monotonicity of the angle means that the curve $\{V(t)| t_i<t<t_{i+1}\}$, leaving at some moment the $\epsilon$-neighborhood of
$V(s)$ in the metric $\angle$ never comes back, or that every ball in $\angle$-metric  with the center $V(s)$ contains only one
connected arc of the considered interval of the curve $V(t)$. Denote by $U(V,\epsilon)$  an  $\epsilon$-neighborhood in the
metric $\angle$ of the point $V$ of $G^p_{p+q}$.

For a given number $\Lambda_5$ there exists some $\omega$ depending on $\Lambda_5$ such that in every $\omega$-ball in the
Riemannian metric $g$ the length of every connected arc of an arbitrary curve with geodesic curvature bounded by $\Lambda_5$ is
bounded by some constant $L$, which depends on $\omega$ and $\Lambda_5$ and has order $2\omega$ as $\omega\to 0$.

Because topologies generated by two metrics  $\angle$ and $g$ coincide, there exists a function $\omega(\epsilon)$ (where
$\omega(\epsilon)\to 0$ as $\epsilon\to 0$) such that every $\epsilon$-ball in the metric $\angle$ is contained in
$\omega(\epsilon)$-ball in the metric $g$. Find some $\epsilon$ such that $\omega(\epsilon)\leq\omega$. Using compactness of
$G^p_{p+q}$ find some finite covering $G^p_{p+q}=\cup_{i} U(V_i,\epsilon/2)$. If some point $V(t')$ belongs to some
$U(V_i,\epsilon/2)$, then because of the triangle inequality the intersection of the considered interval $V(t), t_i<t<t_{i+1}$
with this $U(V_i,\epsilon/2)$ lies in $U(V(t'),\epsilon)$, and by the arguments above has length less than $L$, which obviously
means that the length of the whole interval is bounded by $LN$, where $N$ is the number of all $U(V_i,\epsilon/2)$ in the
finite covering of $G^p_{p+q}$. As was said above, the number of intervals $V(t), t_i<t<t_{i+1}$ with described monotone
behavior of the "angle"-function $\phi_s(t)$ is not bigger than $(p+1)$, so that the length of $V(t),-\infty<t<\infty$ is
bounded by $nLN$, i.e.,
$$
\int_{-\infty}^{\infty}\|\dot V(t)\| dt < \Lambda_6 \tag 2.28
$$
for some constant $\Lambda_6$ depending only on $\Lambda_5$. Lemma~3 and hence, Theorem~3 and Theorem~2 consequently are
proved.
\enddemo

\enddemo

\enddocument

\bye